%%%%%%%%%%%%%%%%%%%%%%%%%%%%%%%%%%%%%%%%%%%%%%%%%%%%%%%%%%%%%%%%%%%%%%%%%%%%%%
%%%%%%%%%%%%%%%%%%%%%%%%%%%%%%%%%%%%%%%%%%%%%%%%%%%%%%%%%%%%%%%%%%%%%%%%%%%%%%
%%%%%%%%%%%%%%%%%%%%%%%%%%%%%%%%%%%%%%%%%%%%%%%%%%%%%%%%%%%%%%%%%%%%%%%%%%%%%%
%%%%%%%%%%%%%%%%%%%%%%%%%%%%%%%%%%%%%%%%%%%%%%%%%%%%%%%%%%%%%%%%%%%%%%%%%%%%%%
%%%%%%%%%%%%%%%%%%%%%%%%%%%%%%%%%%%%%%%%%%%%%%%%%%%%%%%%%%%%%%%%%%%%%%%%%%%%%%
%%%%%%%%%%%%%%%%%%%%%%%%%%%%%%-November-1998-%%%%%%%%%%%%%%%%%%%%%%%%%%%%%%%%%
%%%%%%%%%%%%%%%%%%%%%%%%%%%%%%%%%%%%%%%%%%%%%%%%%%%%%%%%%%%%%%%%%%%%%%%%%%%%%%
%%%%%%%%%%%%%%%%%%%%%%%%%--Beginning-of-ext910.tex--%%%%%%%%%%%%%%%%%%%%%%%%%%
%%%%%%%%%%%%%%%%%%%%%%%%%%%%%%%%%%%%%%%%%%%%%%%%%%%%%%%%%%%%%%%%%%%%%%%%%%%%%%
%%%%%%%%%%%%%%%%%%%%%%%%%%%%%%%%%%%%%%%%%%%%%%%%%%%%%%%%%%%%%%%%%%%%%%%%%%%%%%
%%%%%%%%%%%%%%%%%%%%%%%%%%%%%%%%%%%%%%%%%%%%%%%%%%%%%%%%%%%%%%%%%%%%%%%%%%%%%%

%Preliminaries

\input amstex

\documentstyle{amsppt}

\loadbold

\magnification=\magstep1

\pageheight{9.0truein}
\pagewidth{6.5truein}

\NoBlackBoxes

%%%%%%%%%%%%%%%%%%%%%%%%%%%%%%%%%%%%%%%%%%%%%%%%%%%%%%%%%%%%%%%%%%%%%%%%%%%%%%
%%%%%%%%%%%%%%%%%%%%%%%%%%%%%%%%%%%--MACROS--%%%%%%%%%%%%%%%%%%%%%%%%%%%%%%%%%
%%%%%%%%%%%%%%%%%%%%%%%%%%%%%%%%%%%%%%%%%%%%%%%%%%%%%%%%%%%%%%%%%%%%%%%%%%%%%%
\def\L{{\widetilde L}}
\def\M{{\widetilde M}}
\def\R{\widehat{R}}
\def\grmod{\operatorname{Grmod}}
\def\mod{\operatorname{Mod}}
\def\b{\frak{b}}
\def\Fund{\operatorname{Fund}}
\def\Grfund{\operatorname{GrFund}}
\def\Ext{\operatorname{Ext}}
\def\Extone{\operatorname{Ext}^1}
\def\O{\Cal{O}}

\def\min{\operatorname{min}}

\def\La{L_\alpha}
\def\Lb{L_\beta}
\def\Ma{M_\alpha}
\def\Mb{M_\beta}

\def\plambda{p_\lambda}
\def\Pa{P_\alpha}
\def\Pb{P_\beta}

\def\pa{p_\alpha}
\def\pb{p_\beta}
\def\qa{Q_\alpha}
\def\qb{Q_\beta}

\def\Z{{\Bbb Z}}

\def\ann{\operatorname{ann}}

\def\rank{\operatorname{rank}}
\def\spec{\operatorname{Spec}}
\def\prim{\operatorname{Prim}}

\def\grprim{\operatorname{GrPrim}}
\def\length{\operatorname{length}}

\def\gkdim{\operatorname{GKdim}}
\def\a{\frak{a}}
\def\g{\frak{g}}
\def\go{{\frak{g}}_0}
\def\gi{{\frak{g}}_1}
\def\goi{{\frak{g}}_0 \oplus {\frak{g}}_1}

\def\So{S_0}

\def\Roi{R_0 \oplus R_1}
\def\Ro{R_0}

\def\Xoi{X_0 \oplus X_1}
\def\Xo{X_0}
\def\Xi{X_1}
\def\h{\frak{h}}
\def\ho{\frak{h}_0}

\def\nm{\frak{n}^-}
\def\np{\frak{n}^+}
\def\nmo{\frak{n}^-_0}
\def\npo{\frak{n}^+_0}
\def\U{\boldkey{U}}
\def\V{\boldkey{V}}

\def\Vo{\boldkey{V}_0}

\def\Mo{M_0}
\def\Mi{M_1}
\def\Moi{M_0 \oplus M_1}

\def\ad{\operatorname{ad}}
\def\gr{\operatorname{gr}}
\def\link{{\,\rightsquigarrow\,}}
\def\notlink{\,{{\rightsquigarrow}\kern-0.75em / \kern0.25em}\,}
\def\extlink{{\,\overset \text{ext} \to \rightsquigarrow\,}}
\def\idlink{{\, \thickapprox \,}}

\def\clique{\operatorname{Clique}}
\def\extclique{\operatorname{ext\text{-}Clique}}
%%%%%%%%%%%%%%%%%%%%%%%%%%%%%%%%%%%%%%%%%%%%%%%%%%%%%%%%%%%%%%%%%%%%%%%%%%%%%%
%%%%%%%%%%%%%%%%%%%%%%%%%--REFERENCE-MACROS--%%%%%%%%%%%%%%%%%%%%%%%%%%%%%%%%%
%%%%%%%%%%%%%%%%%%%%%%%%%%%%%%%%%%%%%%%%%%%%%%%%%%%%%%%%%%%%%%%%%%%%%%%%%%%%%%
\def\Bar{{\bf1}}
\def\Bor{{\bf2}}
\def\Bro{{\bf3}}
\def\BroWar{{\bf4}}
\def\CohMon{{\bf5}}
\def\Dix{{\bf6}}
\def\GabJos{{\bf7}}
\def\GooLet{{\bf8}}
\def\GooWar{{\bf9}}
\def\Irv{{\bf10}}
\def\Jan{{\bf11}}
\def\Jat{{\bf12}}
\def\JosSma{{\bf13}}
\def\Kac{{\bf14}}
\def\KraLen{{\bf15}}
\def\LenLet{{\bf16}}
\def\LenWar{{\bf17}}
\def\Letone{{\bf18}}
\def\Lettwo{{\bf19}}
\def\Letthree{{\bf20}}
\def\Letfour{{\bf21}}
\def\Lor{{\bf22}}
\def\McCRob{{\bf23}}
\def\Mon{{\bf24}}
\def\Musone{{\bf25}}
\def\Mustwo{{\bf26}}
\def\PenSer{{\bf27}}
\def\Sch{{\bf28}}
\def\Ser{{\bf29}}
\def\Soe{{\bf30}}
\def\Warone{{\bf31}}
\def\Wartwo{{\bf32}}
%%%%%%%%%%%%%%%%%%%%%%%%%%%%%%%%%%%%%%%%%%%%%%%%%%%%%%%%%%%%%%%%%%%%%%%%%%%%%%
%%%%%%%%%%%%%%%%%%%%%%%%%%--END-OF-MACROS--%%%%%%%%%%%%%%%%%%%%%%%%%%%%%%%%%%%
%%%%%%%%%%%%%%%%%%%%%%%%%%%%%%%%%%%%%%%%%%%%%%%%%%%%%%%%%%%%%%%%%%%%%%%%%%%%%%

\topmatter

%%%%%%%%%%%%%%%%%%%%%%%%%%%%%%%%%%%%%%%%%%%%%%%%%%%%%%%%%%%%%%%%%%%%%%%%%%%%%%
\pretitle{

\hbox{}

\vskip -0.9truein

\noindent {\sevenrm November 1998}

\vskip 0.9truein

         }
%%%%%%%%%%%%%%%%%%%%%%%%%%%%%%%%%%%%%%%%%%%%%%%%%%%%%%%%%%%%%%%%%%%%%%%%%%%%%%

\title Module Extensions Over Classical Lie Superalgebras \endtitle

\rightheadtext{Module Extensions}

\author Edward S. Letzter \endauthor

\address Department of Mathematics, Texas A{\&}M University, College
Station, TX 77843\endaddress

\email letzter\@math.tamu.edu\endemail


\abstract We study certain filtrations of indecomposable injective
modules over classical Lie superalgebras, applying a general approach
for noetherian rings developed by Brown, Jategaonkar, Lenagan, and
Warfield. To indicate the consequences of our analysis, suppose that
$\g$ is a complex classical simple Lie superalgebra and that $E$ is an
indecomposable injective $\g$-module with nonzero (and so necessarily
simple) socle $L$. (Recall that every essential extension of $L$, and
in particular every nonsplit extension of $L$ by a simple module, can
be formed from $\g$-subfactors of $E$.) A direct transposition of the
Lie algebra theory to this setting is impossible. However, we are able
to present a finite upper bound, easily calculated and dependent only
on $\g$, for the number of isomorphism classes of simple highest
weight $\g$-modules appearing as $\g$-subfactors of $E$. \endabstract

\thanks This research was partially supported by grants from the
National Science Foundation. \endthanks

\endtopmatter

\document

\head 1. Introduction \endhead

This paper is about prime ideals, indecomposable injective modules,
and noetherian ring extensions -- building upon \cite{\BroWar, \Jat,
\LenLet, \LenWar, \Lettwo, \Wartwo}. But our work is primarily
motivated by, and applicable to, basic questions arising in the
representation theory of classical Lie superalgebras.

\subhead 1.1 \endsubhead To place our study in context, let $\g$ be a
complex semisimple Lie algebra with Cartan subalgebra $\h$ and Weyl
group $W$. It is well known, within the corresponding category $\O$,
that the $\Ext$ (i.e., $\Extone$) groups of simple modules can be
computed using Kazhdan-Lusztig polynomials \cite{\Bar, \GabJos,
\Irv}. Moreover, calculating $\Ext$ within the category of all
$\g$-modules, the following is elementary: Let $\lambda, \lambda' \in
\h^\ast$. If $L(\lambda)$ and $L(\lambda')$ are simple highest weight
$\g$-modules for which $\Ext \left( L (\lambda), L (\lambda')\right)
\ne 0$, then $\lambda$ and $\lambda'$ occur in the same $W$-orbit in
$\h^\ast$. Of course, there are no nonsplit extensions of finite
dimensional $\g$-modules.

Generalizing further, let $E$ be the injective hull of $L(\lambda)$,
and recall that every essential $\g$-extension of $L(\lambda)$, and in
particular every nonsplit extension of $L(\lambda)$ by a simple
$\g$-module, can be formed from $\g$-subfactors of $E$. It follows
easily, if $L (\lambda')$ is isomorphic to a simple subfactor of $E$,
that $\lambda$ and $\lambda'$ are contained within a single $W$-orbit;
see (5.2). In comparison, if $\g$ is solvable then the injective hull
of a simple finite (and hence one-)dimensional $\g$-module can have
infinitely many pairwise nonisomorphic simple $\g$-subfactors (cf\.
\cite{\BroWar, \S 6}, recalling that winding-automorphism orbits of
prime ideals in this case are either infinite or singletons).

\subhead 1.2 \endsubhead Now, and for the remainder of this section,
let $\g = \goi$ be a complex classical simple Lie superalgebra. (Only
in this section, ``$\g$-module'' will mean ``$\Z_2$-graded
$\g$-module.'') The above theory does not smoothly carry over to this
setting, as is dramatically illustrated by the following example of
Musson \cite{\Mustwo, \S 4}: Over $sl(2,1)$ there exists an infinite
series $L_1,L_2,\ldots$ of finite dimensional simple modules, no two
of which are isomorphic, such that each $\Ext _{sl(2,1)}(L_i,L_{i+1})
\ne 0$. In particular, Kazhdan-Lusztig theory (cf\. \cite{\Ser}) does
not appear to apply. The extension theory for simple highest weight
$\g$-modules (in the sense, e.g., of \cite{\Musone, \PenSer} -- or see
(6.4)) is still largely unknown.

To demonstrate the impact of our study on these issues, let $E$ be an
indecomposable injective $\g$-module with nonzero socle (i.e., $E$ is
the injective hull of a simple $\g$-module), and let $X$ be the set of
simple highest weight $\g$-subfactors (up to isomorphism) of $E$. In
view of Musson's example, and the Lie algebra case, it is natural to
ask whether $X$ can be infinite. We calculate a finite upper bound --
easily determined from $\g$ and not dependent on $E$ -- for the
cardinality of $X$; see (6.7). More precise bounds, concerning
nonsplit extensions of simple highest weight modules, are also given;
see (6.5). Our bounds do not appear to be very sharp, as they
significantly overestimate matters in the case when $\g = \go$.

To more fully describe the results of this paper, in the language of
noetherian ring theory, let $\V$ be the enveloping algebra of $\g$.
The set of {\sl fundamental primes\/} of an injective $\V$-module was
defined in \cite{\LenWar}, following \cite{\BroWar, \Jat}; see
(2.6). Now let $P$ be a prime ideal of $\V$, and let $\Fund (P)$
denote the union of the sets of fundamental primes of indecomposable
injective $\V$-modules with associated prime $P$. Our main theorem
establishes a finite upper bound -- easily calculated from $\g$ and
not dependent on $P$ -- for the cardinality of $\Fund (P)$; see
(6.3). Sharper bounds are found for the number of prime ideals
``linked'' to $P$ by a nonsplit extension of prime submodules; see
(2.3ii) and (6.3).

Our general approach is similar to \cite{\LenLet}, studying ``lying
over''-type properties of fundamental primes in ring extensions. Also,
all of our results for $\g$ depend on the longstanding and well-known
theory for reductive Lie algebras.

\subhead Acknowledgments \endsubhead This work was inspired by
questions raised by Ian Musson --- posed to me in 1991. Much of the
material in \S 2 is folklore, mostly taught to me by Tom Lenagan and
Bob Warfield when I was a graduate student at the University of
Washington.

\head 2. Noetherian Modules \endhead

All of the results discussed in this section should be considered
already known, even in cases when they have not previously appeared
explicitly in the literature; in such instances proofs are included
for completeness. Some of the definitions introduced are nonstandard
but are especially convenient for our later analysis. The reader is
referred to \cite{\GooWar} or \cite{\McCRob} for general background
information on noetherian rings.

\subhead 2.1 \endsubhead Throughout this section, $R$ will denote a
noetherian ring. The set of prime ideals of $R$ will be denoted $\spec
R$ and the set of left primitive ideals of $R$ will be denoted $\prim
R$. The category of $R$-modules will be referred to as $\mod R$, and
$\Ext$ will denote $\Ext^1$. The annihilator in $R$ of a left
$R$-module $M$ will be denoted $\ann _RM$, and the annihilator in $R$
of a right $R$-module $N$ will be denoted $\ann N_R$. Unless otherwise
designated, ``module'' will mean ``left module.''

\subhead 2.2 \endsubhead (Following \cite{\GooWar}.) Let $M$ be a
nonzero  $R$-module.

(i) If $Q$ is a (necessarily prime) ideal of $R$ maximal among
annihilators of nonzero $R$-submodules of $M$, then $\ann Q_M = \{m
\in M \mid Q.m = 0 \}$ is termed an {\sl affiliated submodule\/} of
$M$. A series
%%%%%%%%%%%%%%%%%%%%%%%%%%%%%%%%%%%%%%%%%%%%%%%%%%%%%%%%%%%%%%%%%%%%%
$$0 = M_0 \subset M_1 \subset \cdots \subset M_n = M$$
%%%%%%%%%%%%%%%%%%%%%%%%%%%%%%%%%%%%%%%%%%%%%%%%%%%%%%%%%%%%%%%%%%%%%
of $R$-submodules of $M$ is called an {\sl affiliated series\/} for $M$
if each $M_i/M_{i-1}$ is an affiliated submodule of $M/M_{i-1}$, for
$1 \leq i \leq n$; the prime ideals $\ann _R (M_i/M_{i-1})$ are the
{\sl affiliated primes\/} of the series. More generally, a prime ideal
of $R$ that is an affiliated prime of some affiliated series of $M$ is
referred to as an affiliated prime of $M$.

(ii) If $S$ is a ring and $M$ is an $R$-$S$-bimodule, then the
affiliated series in (i) consists of $R$-$S$-sub-bimodules of $M$; if
$M$ has finite length as an $R$-$S$-bimodule then the set of
affiliated primes of $M$ (viewed as a left $R$-module) coincides with
the set of left annihilators of $R$-$S$-bimodule composition factors
of $M$.

(iii) If a prime ideal $Q$ of $R$ is the annihilator of a submodule of
$M$, then $Q$ is an {\sl annihilator prime\/} of $M$. If $Q = \ann_RM'
= \ann_RM$ for all nonzero $R$-submodules $M'$ of $M$, and $M \ne 0$,
then $M$ is said to be a {\sl prime\/}, or more precisely, a {\sl
$Q$-prime\/} $R$-module. The annihilator of a prime module is
necessarily a prime ideal, and an annihilator of a prime $R$-submodule
of $M$ is termed an {\sl associated prime\/} of $M$. Note that the
annihilator of an affiliated submodule of $M$ is an associated
prime. Also, every nonzero module over a noetherian ring contains a
prime submodule.

(iv) If $M$ is uniform (i.e., does not contain a direct sum of two
nonzero submodules), then it has exactly one associated prime ideal,
which is also the unique ideal of $R$ maximal among annihilators of
nonzero submodules of $M$. If $M$ has finite uniform dimension (i.e.,
contains no direct sum of infinitely many nonzero submodules), then
$M$ can be embedded into a finite direct sum of uniform $R$-module
factors of $M$ such that the associated prime of each uniform factor
is an associated prime of $M$. One way to verify this last statement
is to recall that the injective hull of $M$, which is an essential
extension of $M$ (i.e., $M$ intersects every nonzero $R$-submodule of
its injective hull nontrivially), will be a finite direct sum of
uniform injective $R$-modules; see, for example, \cite{\GooWar,
Chapter 4} for more details.  Also, recall that an injective
$R$-module is indecomposable if and only if it is uniform.

\subhead 2.3 \endsubhead Let $\qa$ and $\qb$ be prime ideals of $R$.

(i) Let $A$ and $B$ be prime noetherian rings, and suppose that $C$ is
a nonzero $A$-$B$-bimodule finitely generated and torsion free on each
side; we say that $C$ is a {\sl bond\/} from $A$ to $B$. (A bimodule
finitely generated on each side, over prime noetherian rings, is
torsion free on each side if and only if it is faithful and prime on
each side; see, e.g., \cite{\Jat, 5.1.1}.) It follows, for example,
from \cite{\GooWar, 7.16} that if $A$ is right primitive then $B$ is
right primitive, and if $B$ is left primitive then $A$ is left
primitive. Now suppose there are ideals $I \subset J$ of $R$ such that
$\ann _R (I/J) = \qa$, such that $\ann (I/J)_R = \qb$, and such that
$I/J$ is a bond from $R/\qa$ to $R/\qb$. Then we say that $I/J$ is an
{\sl ideal link\/} from $\qa$ to $\qb$ and write $\qa \idlink
\qb$. (But ``$\idlink$'' is not necessarily an equivalence relation.)
If $\qa\qb \subseteq I \subset J = \qa \cap \qb$, then we say that
there is a {\sl link\/} from $\qa$ to $\qb$ and write $\qa \link
\qb$. Occasionally, for emphasis, we will refer to links as {\sl
direct links\/}. If there exist prime ideals $Q_1,\ldots,Q_m$ and a
sequence of links $\qa \link Q_1 \link \cdots \link Q_m \link \qb$,
then we write $\qa \link \cdots \link \qb$. (When writing $\qa \link
\cdots \link \qb$, we include the possibility that $\qa = \qb$ and no
link exists from $\qa$ to itself.)  If $Q$ is a prime ideal of $R$,
then the set of all prime ideals $Q'$ such that either $Q' \link
\cdots \link Q$ or $Q \link \cdots \link Q'$ is termed the {\sl
clique\/} of $Q$ and denoted $\clique (Q)$.

(ii) We will say that a finitely generated uniform  $R$-module $M$
is an {\sl external link\/} from $\qa$ to $\qb$, and write $\qa
\extlink \qb$, provided that the following condition holds: There
exists a short exact sequence
%%%%%%%%%%%%%%%%%%%%%%%%%%%%%%%%%%%%%%%%%%%%%%%%%%%%%%%%%%%%%%%%%%%%%%%
$$0 \rightarrow \Ma \rightarrow M \rightarrow \Mb \rightarrow
0$$
%%%%%%%%%%%%%%%%%%%%%%%%%%%%%%%%%%%%%%%%%%%%%%%%%%%%%%%%%%%%%%%%%%%%%%%
of uniform  $R$-modules such that $0 \subset \Ma \subset M$
is an affiliated series, such that $\qa = \ann _R\Ma$, such that
$\qb = \ann _R\Mb$, and such that $\ann_RM' = \ann_RM$ for all
$R$-submodules $M'$ of $M$ not contained in $\Ma$. When there are
prime ideals $Q_1,\ldots,Q_m$ and a sequence of external links $\qa
\extlink Q_1 \extlink \cdots \extlink Q_m \extlink \qb$, we write $\qa
\extlink \cdots \extlink \qb$. (As before, $\qa \extlink \cdots
\extlink \qa$ whether or not $\qa \extlink \qa$.) If $Q$ is a prime
ideal of $R$, the set of all prime ideals $Q'$ such that either $Q
\extlink \cdots \extlink Q'$ or $Q' \extlink \cdots \extlink Q$ is
termed the {\sl external clique\/} of $Q$ and denoted $\extclique
(Q)$. 

(iii) (Cf. \cite{Jat, 6.1.3}.) Let $M$ be an external link from $\qa$
to $\qb$, and set $I = \ann _RM$. Then Jategaonkar's Main Lemma (as
stated, e.g., in \cite{\GooWar, 11.1}) asserts that one of the
following two cases must occur: (1) $\qa \cap \qb/I$ is a link from
$\qa$ to $\qb$, or (2) $\qb = I \subsetneq \qa$. 

It follows from \cite{\Bro, 2.3ii} (cf\. \cite{\GooWar, 11.2}) that if
$\qa \link \qb$ then $\qa \extlink \qb$.

(iv) The condition in (ii) can be weakened somewhat. Let $N$ be a
uniform  $R$-module, and suppose that $\qa$ is the unique
associated prime ideal of $N$. As noted previously, $\qa$ is the
unique maximal annihilator in $R$ of nonzero submodules of $N$. Set
$N_\alpha = \ann (\qa)_N$. Suppose further that $0 \subset N_\alpha
\subset N$ is an affiliated series, and that $\qb = \ann
N/N_\alpha$. Setting $I$ equal to the ideal of $R$ maximal among
annihilators of submodules not containing $N_\alpha$, we see that $M =
\ann I_N$ is an external link from $\qa$ to $\qb$.

(v) If for all external links between prime ideals of $R$, only case
(1) in (iii) holds, then $R$ satisfies the {\sl (left) strong second
layer condition}.

(vi) Suppose that $R$ satisfies the strong second layer condition and
that $\qa \idlink \qb$. Then $\qa \link \cdots \link \qb$; see
\cite{\Jat, 8.2.4}.

\subhead\nofrills\endsubhead

Although the following lemma is probably well known, we do not know of
a suitable reference. The proof is adapted from \cite{\Jat, 7.1.2},
which applies to the case when $R$ satisfies the strong second layer
condition.

\proclaim{2.4 Lemma} Let $M$ be a finitely generated nonzero
$R$-module, and let $X$ denote the union of the external cliques of
the associated primes of $M$. Then there exist $Q_1,\ldots,Q_t \in X$
such that $Q_i$ is the annihilator of a prime $R$-module subfactor of
$M$, for $1 \leq i \leq t$, and such that $Q_1\cdots Q_tM = 0$. If $M$
has finite length then $\{ Q_1,\ldots,Q_t \}$ can be chosen to be the
set of annihilators of composition factors of $M$. \endproclaim

\demo{Proof} By noetherian induction, we may assume that the
conclusion of the lemma holds for all of the proper $R$-module factors
of $M$. Now suppose that $M$ is not uniform. As noted in (2.2iv), $M$
embeds into a finite direct sum of proper uniform $R$-module factors
of $M$, each of whose associated primes is an associated prime of $M$;
the lemma follows in this case. Therefore, we may assume that $M$ is
uniform. Let $Q$ denote the unique associated prime of $M$, and let $N
= \ann_MQ$. If $M$ has finite length then $Q$ is the annihilator of a
simple submodule of $N$. If $Q'$ is any associated prime of $M/N$,
then by (2.3iv) we can find a submodule of (the uniform module) $M$
providing an external link from $Q$ to $Q'$. Therefore, since the
conclusion of the lemma holds for $M/N$, the desired conclusion also
holds for $M$.  \qed\enddemo

The following is certainly well known, but as above we do not know of
an appropriate reference.

\proclaim{2.5 Lemma} Suppose that $\qa$ and $\qb$ are prime ideals of
$R$ such that $\qa \extlink \qb$. Then $\qa \cap Z = \qb \cap Z$,
where $Z$ is the center of $R$. \endproclaim

\demo{Proof} Let $M$ be an external link from $\qa$ to $\qb$, as in
(2.3ii). When $\qa \link \qb$, the conclusion of the lemma follows
immediately. Therefore, we may assume without loss of generality that
$\qa \supsetneq \qb = \ann_RM$, as in (2.3iii), and we may further
assume without loss that $\qb = 0$. (Note, however, that the image in
$R/\qb$ of the original $Z$ may be smaller than the center of
$R/\qb$.)  Now let $I = \qa \cap Z$, and set $J = RI$. By
\cite{\GooWar, 11.13}, for example, $J$ is an AR-ideal (see, e.g.,
\cite{\GooWar, p\. 190}). Thus $J^nM = 0$ for some positive integer
$n$, by \cite{\GooWar, 11.11b}, since $J$ annihilates a nonzero
submodule of $M$, and since $M$ is uniform. Hence $I^nM = 0$. But $M$
is a faithful $Z$-module, and $Z$ is prime. Consequently, $I = J = 0$,
and the lemma follows. \qed\enddemo

\subhead 2.6 \endsubhead We now briefly review the construction and
basic properties of the fundamental series of an $R$-module, following
\cite{\LenWar} (cf\. \cite{\BroWar; \Jat, Chapter 9}). To start,
assume that $R$ is an algebra of finite Gelfand-Kirillov dimension
(GK-dimension) over a field $k$ (cf\. \cite{\KraLen} or \cite{\McCRob}
for definitions and background), that the GK-dimension of every
$R$-module with finite GK-dimension is an integer, and that if $0
\rightarrow L \rightarrow M \rightarrow N \rightarrow 0$ is an exact
sequence of $R$-modules then the GK-dimension of $M$ is equal to the
maximum of the GK-dimensions of $L$ and $N$. (The preceding
assumptions are summarized by saying that $R$ has {\sl exact integer
GK-dimension\/}.)

(i) Let $M$ be a nonzero $R$-module. Say that $M$ is {\sl
($\alpha$-)annihilator homogeneous\/} for some integer $\alpha$ if
$\gkdim (R/\ann_RN) = \alpha$ for all nonzero finitely generated
submodules $N$ of $M$. If $X$ is a set of semiprime ideals of $R$,
then the intersection of the ideals in any nonempty subset of $X$ is
referred to as an {\sl $X$-semiprime\/} ideal of $R$.

(ii) Assume that $M$ is an annihilator homogeneous nonzero 
$R$-module, and set $F_0(M) = 0$. For $j = 1,2,\ldots$, set $X_j(M)$
equal to the set of associated primes of $M/F_{j-1}(M)$, and set
%%%%%%%%%%%%%%%%%%%%%%%%%%%%%%%%%%%%%%%%%%%%%%%%%%%%%%%%%%%%%%%%%%%%%%%%%%
$$F_j(M) = \left\{ m \in M \mid \text{$Im \in F_{j-1}(M)$ for some
$X_j(M)$-semiprime ideal $I$ of $R$}\right\} .$$
%%%%%%%%%%%%%%%%%%%%%%%%%%%%%%%%%%%%%%%%%%%%%%%%%%%%%%%%%%%%%%%%%%%%%%%%%%
The series $F_0(M) \subset F_1(M) \subset \cdots$ is the {\sl
fundamental series\/} for $M$. In \cite{\LenWar, 3.2} it is shown that
$M = \cup _j F_j(M)$. The set of {\sl fundamental primes\/} of $M$,
denoted $\Fund (M)$, is defined as $\cup_j X_j(M)$.

(iii) Now let $M$ be an arbitrary (i.e., not necessarily annihilator
homogeneous) nonzero  $R$-module, and let $\alpha_1 < \alpha_2 <
\cdots < \alpha_n$ be the distinct integers arising as the
GK-dimensions of $R/ \ann _RN$, for finitely generated submodules $N$
of $M$. Set $G_0(M) = 0$, and for $1 \leq i \leq n$ set $G_i(M)$ equal
to
%%%%%%%%%%%%%%%%%%%%%%%%%%%%%%%%%%%%%%%%%%%%%%%%%%%%%%%%%%%%%%%%%%%%%%%%%%%
$$\left\{ m \in M \mid \text{$m \in G_{i-1}(M)$ or
$\gkdim\left(R/\ann_R\left[Rm+G_{i-1}(M)/G_{i-1}(M)\right]\right) =
\alpha_i$} \right\} .$$
%%%%%%%%%%%%%%%%%%%%%%%%%%%%%%%%%%%%%%%%%%%%%%%%%%%%%%%%%%%%%%%%%%%%%%%%%%%
In \cite{\LenWar, 4.2} it is proved that
%%%%%%%%%%%%%%%%%%%%%%%%%%%%%%%%%%%%%%%%%%%%%%%%%%%%%%%%%%%%%%%%%%%%%%%%%%%
$$G_0(M) \subsetneq G_1(M) \subsetneq \cdots \subsetneq G_n(M) = M ,$$
%%%%%%%%%%%%%%%%%%%%%%%%%%%%%%%%%%%%%%%%%%%%%%%%%%%%%%%%%%%%%%%%%%%%%%%%%%%
and that each $G_i(M)/G_{i-1}(M)$ is an $\alpha_i$-annihilator
homogeneous $R$-module. Set
%%%%%%%%%%%%%%%%%%%%%%%%%%%%%%%%%%%%%%%%%%%%%%%%%%%%%%%%%%%%%%%%%%%%%%%%%%%
$$\Fund _i(M) =\Fund (G_i(M)/G_{i-1}(M)),$$
%%%%%%%%%%%%%%%%%%%%%%%%%%%%%%%%%%%%%%%%%%%%%%%%%%%%%%%%%%%%%%%%%%%%%%%%%%%
and define the {\sl fundamental primes\/} of $M$, denoted $\Fund (M)$,
to be $\cup _i \Fund_i(M)$.

(iv) Continue to let $M$ be an arbitrary nonzero  $R$-module, and
let $Q$ be an annihilator prime of $M$. In \cite{\LenWar, 4.4} it is
proved that there exist submodules $H$ and $J$ of $M$ such that $QH =
0$ and such that $H/J$ is $Q$-prime. Choosing $h \in H \setminus J$,
we see that $Q = \ann _RR.h$. Consequently, the set of annihilator
primes of $M$ is equal to the set of prime annihilators of finitely
generated submodules of $M$.

\subhead 2.7 \endsubhead (i) Let $A$ and $B$ denote algebras of finite
GK-dimension, and let $C$ denote an $A$-$B$-bimodule finitely
generated on each side. Then $\gkdim (A/\ann _AC) = \gkdim (_AC) =
\gkdim (C_B) = \gkdim (B/\ann C_B)$; see, for example, \cite{\KraLen,
5.3}. Consequently, if $A$ is a subalgebra of $B$, and if $B$ is
finitely generated on either the right or left as an $A$-module, then
$\gkdim (A) = \gkdim(B)$.

(ii) Suppose that $R$ has finite GK-dimension. Let $\qa$ and $\qb$ be
prime ideals of $R$. If $\qa \subseteq \qb$ then $\qa = \qb$ if and
only if $\gkdim (R/\qa) = \gkdim (R/\qb)$; see \cite{\KraLen, 3.16}.

\subhead 2.8 \endsubhead Assume that $R$ has exact integer
GK-dimension. Let $E$ be an injective $R$-module, and let $\qb$ be a
prime ideal of $R$. It is proved in \cite{\LenWar, 5.4} (noting
(2.7i)) that $\qb$ is a fundamental prime of $E$ if and only if there
exists an annihilator prime $\qa$ of $E$, and an ideal link $I/J$ from
$\qa$ to $\qb$, such that $\gkdim (R/\qa) = \gkdim (R/J) = \gkdim
(R/\qb)$.

\subhead 2.9 \endsubhead Suppose that $R$ has finite GK-dimension. Let
$\qa$ and $\qb$ be prime ideals of $R$. 

(i) In view of (2.8), we will say that there is a {\sl strong\/}
ideal link from $\qa$ to $\qb$ if there is an ideal link $I/J$ such
that $\gkdim (R/\qa) = \gkdim (R/J) = \gkdim (R/\qb)$.

(ii) Suppose that there is an ideal link $I/J$ from $\qa$ to $\qb$ and
that $\gkdim (R/\qa) = \gkdim (R/\qb) = \gkdim (R)$. Since $J$ is
contained in both $\qa$ and $\qb$, we see that $\gkdim (R) \geq \gkdim
(R/J) \geq \gkdim (R/\qa) = \gkdim (R)$. Therefore, $I/J$ is a strong
ideal link.

(iii) Now assume that $R$ has exact integer GK-dimension. It follows
from the remarks after \cite{\LenWar, 5.4}, when there exists a strong
ideal link from $\qa$ to $\qb$, that $\qa \link \cdots \link
\qb$. Hence if $E$ is an injective  $R$-module, and $\qb$ is a
fundamental prime of $E$, then there exists an annihilator prime $\qa$
of $E$ such that $\qa \link \cdots \link \qb$.

\subhead 2.10 \endsubhead Assume that $R$ has exact integer
GK-dimension, and suppose that $(\qa \cap \qb)/J$ is a link from $\qa$
to $\qb$. It follows from (2.7i) and, for example, \cite{\KraLen, 5.7}
that $\gkdim (R/\qa) = \gkdim (R/J) = \gkdim (R/\qb)$. We see in this
case that a direct link is a strong ideal link.

\proclaim{2.11 Lemma} Assume that $R$ has exact integer GK-dimension,
that $M$ is a uniform  $R$-module, that $\qa$ is the unique
associated prime of $M$, and that $\qb$ is an annihilator prime of
$M$. Then $\qa \extlink \cdots \extlink \qb$. \endproclaim

\demo{Proof} By (2.6iv), $\qb$ is the annihilator of a finitely
generated submodule of $M$; this submodule is uniform and its unique
associated prime is $\qa$. Therefore, we can reduce to the case where
$\qb = 0$ and where $M$ is a finitely generated faithful $R$-module.
However, by (2.4), there exist prime ideals $Q_1,\ldots,Q_t$ in the
external clique of $\qa$ such that $Q_1\cdots Q_tM = 0$. Thus
$Q_1\cdots Q_t = 0$, and so $Q_i = 0$ for some $1 \leq i \leq t$. The
lemma follows. \qed\enddemo

\subhead 2.12 \endsubhead Continue to assume that $R$ has exact
integer GK-dimension, and let $Q$ be a prime ideal of $R$. 

(i) Define $\Fund (Q)$, the set of {\sl fundamental primes\/} of $Q$,
to be the union of the sets of fundamental primes of indecomposable
(equivalently, uniform) injective $R$-modules for which $Q$ is the
unique associated prime. It follows from (2.9iii) and (2.11) that
$\Fund (Q) \subseteq \extclique (Q)$.

(ii) If $M$ is a finite length uniform  $R$-module with associated
primitive ideal $Q$, then it is not hard to verify that the
annihilators of the composition factors of $M$ are all contained in
$\Fund(Q)$.

(iii) Suppose that $M$ is a uniform  $R$-module with associated
prime $Q$, that $L$ is an arbitrary simple $R$-module subfactor of
$M$, and that $\ann _RL = \qb$. We may embed $M$ into a uniform
injective $R$-module $E$, and the filtration of $E$ described in (2.6)
is exhaustive. Hence there exists $\qa \in \Fund(Q) \subseteq
\extclique (Q)$ such that $\qa \subseteq \qb$.

\head 3. Noetherian Superalgebras \endhead

We now briefly review some generalities, relevant to the study of Lie
superaglebras, from the theory of $\Z_2$-graded rings.

\subhead 3.1 \endsubhead (i) Throughout this section, $k$ will denote
a field not of characteristic two, and $R$ will denote a noetherian
$k$-algebra. We will further assume, throughout, that $R = \Roi$ is
$\Z _2$-graded (i.e., $R$ is an associative superalgebra). 

(i) Henceforth, we will use the term {\sl graded\/} to mean
``$\Z_2$-graded'' and ``homogeneous with respect to the $\Z _2
$-grading.'' If $X = \Xoi$ is a ($\Z_2$-)graded $k$-vector space, then
the elements of $\Xo$ will be called {\sl even\/} and the elements of
$\Xi$ will be called {\sl odd\/}. Degree-preserving linear maps
between graded $k$-vector spaces will be termed {\sl
graded\/}. However, a $k$-vector space should be assumed graded only
when explicitly specified. {\sl Ungraded\/} will mean ``not
necessarily graded.''

(ii) Let $X = \Xoi$ be a graded $k$-vector space. There is an
automorphism $\sigma$ of $X$, obviously of order two, mapping each
even element to itself and each odd element to its additive
inverse. Also, $\Xo$ is identical to $X^\sigma$, the subspace of
$\sigma$-invariants. Now let $Y$ be an ungraded subspace of $X$. The
unique maximum graded subspace of $X$ contained in $Y$ will be termed
$\gr Y$, and it is easy to see that $Y$ is graded if and only if
$\sigma (Y) = Y$. Therefore, $\gr Y = Y \cap \sigma (Y)$.

(iii) Applied to $R$, $\sigma$ is a $k$-algebra
automorphism. Therefore, $\Ro$ is noetherian and $R$ is finitely
generated on each side as an $\Ro$-module \cite{\Mon, 1.12,
5.9}. Applied to a graded $R$-module, $\sigma$ permutes the ungraded
submodules.

(iv) The category of graded left $R$-modules -- with graded morphisms
-- will be denoted $\grmod R$.

\subhead 3.2 \endsubhead (i) We will say that a graded $R$-module $M$
is {\sl graded-simple\/} provided the only nonzero graded submodule of
$M$ is itself (i.e., $M$ is a simple object in $\grmod R$), a
condition that occurs if and only if $\Mo$ and $\Mi$ are simple
$\Ro$-modules.  Consequently, a graded $R$-module $N$ has finite
length as an ungraded $R$-module if and and only if $N$ has finite
length in $\grmod R$.

(ii) The annihilators of  graded-simple $R$-modules are the {\sl
(left) graded-primitive\/} ideals, and the set of graded-primitive
ideals of $R$ will be denoted $\grprim R$.

(iii) Assume that $M$ is a graded-simple  $R$-module with
annihilator $p$. The $\Ro$-module structure of $M$ ensures that its
length as an ungraded $R$-module is no greater than two. Moreover, if
$N$ is an ungraded simple $R$-submodule of $M$, and $\ann _RN = P$,
then $M = N + \sigma(N)$, and $p = \gr (P) = P\cap \sigma (P)$.

(iv) Let $L$ be a maximal ungraded left ideal of $R$. Set $P = \ann
_R(R/L)$, and set $M = R/L \cap \sigma(L)$. Then $M \subseteq
(R/L)\oplus (R/\sigma (L))$ has length no greater than two, as either
an ungraded or graded $R$-module, and $p = \gr (P) = \ann _RM$. It is
now straightforward to check that $p$ is graded-primitive.

(v) We see from (iii) and (iv) that a graded ideal $p$ of $R$ is
graded-primitive if and only if there exists a primitive ideal $P$ of
$R$ for which $p = \gr (P) = P \cap \sigma (P)$.

\subhead 3.3 \endsubhead (i) A graded ideal of $R$ is {\sl
graded-prime\/} when it contains no product of strictly larger graded
ideals, and a graded ideal is {\sl graded-semiprime\/} when it is the
intersection of graded-prime ideals. Also, $R$ will be termed {\sl
graded-prime\/} (resp\. {\sl graded-semiprime\/}) when the zero ideal
of $R$ is graded-prime (resp\. graded-semiprime). It is proved as
follows that an ideal $p$ of $R$ is graded-prime if and only if there
exists a prime ideal $P$ of $R$ such that $p = \gr (P) = P \cap \sigma
(P)$: First, when $P$ is a prime ideal of $R$, it is straightforward
to check that $\gr (P)$ is graded-prime. Conversely, let $p$ be a
graded prime ideal of $R$. Since $R$ is noetherian there exist prime
ideals $P_1,\ldots,P_t$ of $R$, each containing $p$, such that
$P_1\cdots P_t \subseteq p$. But $\gr (P_1)\cdots \gr (P_t) \subseteq
p$, and each of the $\gr (P_i)$ contains $p$. Thus $p = \gr (P_i)$ for
some $1 \leq i \leq t$.

(ii) If $P$ is a prime ideal of $R$, then $P \cap \Ro = \sigma(P) \cap
\Ro = \gr (P) \cap \Ro = P^\sigma$ is a semiprime ideal of $\Ro$; see,
for example, \cite{\Mon, 1.5}. 

(iii) Suppose that $R$ is graded-semiprime. Then $\Ro$ is semiprime,
by (ii). Let $X$ denote the Ore set of regular elements of $\Ro$. By
\cite{\Mon, 5.3}, $X$ is an Ore set of regular elements of $R$, and
the left Ore localization $X^{-1}R$ is equal to the Goldie quotient
ring of $R$. Because $X$ consists only of even elements of $R$, the
grading on $R$ extends naturally to a grading of $X^{-1}R$.

Now let $M$ be a graded  $R$-module, and let $Y$ denote the Ore
set of regular elements of $R$.  Because $X^{-1}R$ is the Goldie
quotient ring of $R$, a submodule $N$ of $M$ is $X$-torsion (i.e.,
torsion with respect to $X$) if and only if $N$ is torsion (i.e.,
$Y$-torsion). Since $X$ consists only of even elements, the torsion
submodule of $M$ is therefore graded. Also, $M$ is torsion free if and
only if $M$ is $X$-torsion free.

(iv) Continue to assume that $R$ is graded-semiprime. Let $M=\Moi$ be
a graded $R$-$R$-bimodule (i.e., $R_iM_j$ and $M_jR_i$ are both
contained in $M_{i+j}$). By (iii), the left torsion submodule $T$ of
$M$ is a graded bimodule, and by \cite{\Jat, 5.1.1}, $T$ is the
largest left unfaithful graded $R$-$R$-sub-bimodule of $M$. We may
therefore define a {\sl graded-bond\/} between graded-prime rings $A$
and $B$ to be a graded bimodule finitely generated and torsion free on
each side. {\sl Graded-ideal-links\/} and {\sl
graded(-direct-)links\/} can be similarly defined.

(v) Suppose that $R$ is graded-prime. Then its Goldie quotient ring is
a graded-simple (i.e., containing no proper graded ideals) Artinian
ring. Now suppose that $A$ is a graded-simple artinian ($\Z_2$-)graded
ring. By (i), $A$ is either simple artinian or equal to a direct
product of two isomorphic simple artinian rings. In the first case,
$A$ is an $m\times m$ matrix ring over a division ring and we say that
$m$ is the {\sl graded-rank\/} of $A$. In the second case, $A$ is a
direct product of two $n \times n$ matrix rings over division rings,
and we say that $n$ is the {\sl graded-rank\/} of $A$. The {\sl
graded-Goldie-rank\/} of $R$ is the graded-rank of its Goldie quotient
ring.

\subhead 3.4 \endsubhead Assume that $R$ has exact integer
GK-dimension -- see (2.6) -- and that $M$ is a graded  $R$-module.

(i) Let $I$ be an ideal of $R$. Since $R$ has exact integer
GK-dimension, $\gkdim (R/I) = \gkdim (R/\sigma (I)) = \gkdim [R/(I\cap
\sigma(I))] = \gkdim (R/ \gr (I))$, by \cite{\KraLen, 5.7}. If $m
\in M$ and $I = \ann _R R m$, then $\sigma (I) = \ann _R R \sigma (m)$.

(ii) Let $\alpha$ be the least positive integer occurring among
$\gkdim _R (N)$ for finitely generated nonzero $R$-submodules $N$ of
$M$. It follows from (i) that $G(M) = \{ m \in M \mid \gkdim (R / \ann
_R R m) = \alpha \}$ is a $\sigma$-stable $R$-submodule of $M$, and so
$G(M)$ is a graded submodule of $M$. In particular, $G(M)$ is equal to
the sum of homogeneous $m \in M$ such that $\gkdim (R / \ann _R R m) =
\alpha$.

(iii) Let $\alpha_1 < \alpha_2 < \cdots < \alpha_n$ be the integers
occurring among $\gkdim (R /\ann _R N)$ for finitely generated
nonzero submodules $N$ of $M$. By (ii), each of the submodules in the
series $0 = G_0(M) \subset G_1(M) \subset \cdots \subset G_n(M)= M$,
as defined in (2.6iii), is graded. 

(iv) Suppose that $\gkdim (R /\ann _R L) = \alpha$ for all finitely
generated nonzero graded $R$-submodules $L$ of $M$; we will say in
this case that $M$ is {\sl
($\alpha$-)graded-annihilator-homogenous\/}. Let $N$ be a finitely
generated nonzero ungraded submodule of $M$, and suppose that $\ann _R
N = I$. Then $\alpha = \gkdim [R/ \ann (N + \sigma (N)] = \gkdim [R /
\gr (I)] = \gkdim (R/I)$, by (i). Hence $M$ is $\alpha$-annihilator
homogeneous (see (2.6i)).

(v) We see from (iv) that a graded $R$-module is
$\alpha$-graded-annihilator homogeneous if and only if it is
$\alpha$-annihilator homogeneous. For $1 \leq i \leq n$, the graded
subfactors $G_i(M)/G_{i-1}(M)$ of (iii) are $\alpha
_i$-graded-annihilator-homogeneous, by (2.6iii).

\subhead 3.5 \endsubhead The following duality principle of Cohen and
Montgomery \cite{\CohMon} allows us to apply ungraded noetherian ring
theory to $\grmod R$. Set $\R = R\# k [\Z_2 ]^*$ as in \cite{\CohMon,
\S 1}. When viewed as either a left or right module over its
$k$-subalgebra $R$, $\R$ is free of rank two. Hence $\R$ is
noetherian, and if $R$ has finite GK-dimension then so does $\R$.
Furthermore, if $R$ has exact integer GK-dimension then so does $\R$,
by \cite{\Lor, 1.6}.

In \cite{\CohMon, 2.2} it is proved that the categories $\grmod R$ and
$\mod \R$ are isomorphic.

\subhead 3.6 \endsubhead (i) In view of (3.2--5), we now claim that
the results mentioned in \S 2, and their prerequisites, are valid if
the modules and bimodules involved are replaced with the corresponding
graded objects. (The terms ``dimension over $k$'' and ``GK-dimension''
retain their usual ungraded definitions even when applied to graded
modules.)  A formal proof of this assertion will be omitted.

(ii) We will continue to use the prefix {\sl graded-\/} to denote the
graded analogue of an ungraded term.

(iii) The existence of a graded-external-link from a graded-prime
ideal $\pa$ to a graded-prime ideal $\pb$ will be denoted $\pa
\extlink \pb$, and a graded-link from $\pa$ to $\pb$ will be denoted
$\pa \link \pb$. When $R$ has exact integer GK-dimension, $\Grfund
(M)$ will refer to the fundamental graded-primes of a graded
$R$-module $M$, and if $p$ is a graded-prime ideal of $R$ we will use
$\Grfund (p)$ to denote the union of the sets of the fundamental
graded-primes of graded-uniform graded-injective $R$-modules.

(iv) Note that the graded-injective hull of a graded  $R$-module
$M$ will be a maximal graded-essential extension of $M$.

\head 4. Extensions of Rings \endhead

This section is devoted to developing some preparatory results
concerning ring extensions finitely generated on one side.

\subhead 4.1 \endsubhead Throughout this section $S$ will denote a
noetherian algebra of finite GK-dimension over a field $k$, and $R$
will denote a noetherian subalgebra of $S$ such that $S = s_1R +
\cdots + s_mR$ for some $s_1,\ldots,s_m \in S$. Further suppose,
throughout, that $P$ is a prime ideal of $S$ and that $Q$ is a prime
ideal of $R$.

(i) Suppose that $Q$ is minimal over $P\cap R$; we say that $P$ {\sl
lies over\/} $Q$ and that $Q$ {\sl lies under\/} $P$. Since $R$ is
noetherian, at least one and at most finitely many prime ideals of $R$
lie under $P$. It follows from \cite{\Letone, 1.1}, for example, that
there is a bond from $S/P$ to $R/Q$. Consequently, by (2.7i), $\gkdim
(R/Q) = \gkdim (S/P)$.  Of course, it is easy to see that $P$ has
finite codimension (in $S$) if and only if $Q$ has finite codimension
(in $R$). From (2.3i) it follows that if $P$ is right primitive then
$Q$ is right primitive, and if $Q$ is left primitive then $P$ is left
primitive.

(ii) If $P$ lies over $Q$, and if $Q$ is a left annihilator prime of
$S/P$, then we say that $P$ {\sl lies directly over\/} $Q$ and that
$Q$ {\sl lies directly under\/} $P$. (Cf\. \cite{\GooLet, 5.1}.) There
exists at least one prime ideal of $R$ lying directly under $P$ (see,
e.g., \cite{\Letfour, 2.6i}).

\subhead 4.2 \endsubhead Let $F$ denote the Goldie quotient ring of
$R/Q$. Because $S$ is finitely generated as a right $R$-module,
$S\otimes _RF$ has finite length as a right $F$-module. It is proved
in \cite{\Letfour, 2.4} that $P$ lies directly over $Q$ if and only if
$P$ is the left annihilator in $S$ of an irreducible $S$-$F$-bimodule
factor of $S\otimes _RF$.

The proof of the following relies on the Joseph-Small Additivity
Principle (cf\. \cite{\Bor, \JosSma, \Warone}).

\proclaim{4.3 Proposition} {\rm (i)} $P$ lies directly over at most $m$
prime ideals of $R$. {\rm (ii)} $Q$ lies directly under at most $m$
prime ideals of $S$. \endproclaim

\demo{Proof} (i) Let $Q_1,\ldots,Q_s$ be prime ideals of $R$ lying
directly under $P$. For $1 \leq i \leq s$, let $F_i$ denote the Goldie
quotient ring of $R/Q_i$. It follows from (4.2) that $P$ is the left
annihilator of an irreducible $S$-$F_i$-bimodule factor of $S\otimes
_R F_i$. By \cite{\GooWar, 7.23}, a version of the Additivity
Principle,
%%%%%%%%%%%%%%%%%%%%%%%%%%%%%%%%%%%%%%%%%%%%%%%%%%%%%%%%%%%%%%%%%%%%%%%
$$\rank (S/P) \leq \length (S \otimes _RF_i)_{F_i} \leq
m\cdot\rank(R/Q_i),$$
%%%%%%%%%%%%%%%%%%%%%%%%%%%%%%%%%%%%%%%%%%%%%%%%%%%%%%%%%%%%%%%%%%%%%%%
where $\rank()$ denotes Goldie rank and $\length()_{F_i}$ denotes
composition length as a right $F_i$-module. However, it now follows
from the standard version of the Additivity Principle (see, e.g.,
\cite{\McCRob, 4.5.4}) that
%%%%%%%%%%%%%%%%%%%%%%%%%%%%%%%%%%%%%%%%%%%%%%%%%%%%%%%%%%%%%%%%%%%%%%%
$$\rank(R/Q_1) + \cdots + \rank(R/Q_s) \leq \rank (S/P)) \leq m\cdot
\min \left\{\rank(R/Q_1),\ldots,\rank(R/Q_s)\right\}.$$
%%%%%%%%%%%%%%%%%%%%%%%%%%%%%%%%%%%%%%%%%%%%%%%%%%%%%%%%%%%%%%%%%%%%%%%
Hence $s \leq m$.

(ii) Suppose that $P_1,\ldots,P_t$ are prime ideals of $S$ lying
directly over $Q$. Another application of the Additivity Principle
(this time see, e.g., \cite{\GooWar, 7.26}) ensures that $\rank (R/Q)
\leq \rank (S/P_j)$, for $1 \leq j \leq t$. Let $F$ denote the Goldie
quotient ring of $R/Q$; it follows from (4.2) that $P_j$ is the
annihilator of an irreducible $S$-$F$-bimodule factor of $S\otimes
_RF$. Therefore, again using \cite{\GooWar, 7.23}, we see that
%%%%%%%%%%%%%%%%%%%%%%%%%%%%%%%%%%%%%%%%%%%%%%%%%%%%%%%%%%%%%%%%%%%%%%%
$$\rank (S/P_1) + \cdots + \rank (S/P_t) \leq m\cdot \min\left\{ \rank
(S/P_1), \ldots, \rank(S/P_t)\right\},$$
%%%%%%%%%%%%%%%%%%%%%%%%%%%%%%%%%%%%%%%%%%%%%%%%%%%%%%%%%%%%%%%%%%%%%%%
and so $t \leq m$. \qed\enddemo

The following is analogous to \cite{\Letthree,
3.4}.

\proclaim{4.4 Proposition} Suppose that $P$ lies over $Q$. Then there
exists a prime ideal $Q'$ of $R$ such that $Q'$ lies directly under
$P$ and such that there is a strong ideal link from $Q'$ to
$Q$. \endproclaim

\demo{Proof} We may assume, without loss of generality, that $P =
0$. We first find a prime ideal $Q'$ of $R$ such that $Q'$ is an
annihilator prime of $S$ and such that there is an ideal link from
$Q'$ to $Q$. Let $X$ denote the set of elements of $R$ regular modulo
its nilradical. Because $S$ has finite GK-dimension and is finitely
generated as a right $R$-module, it follows from \cite{\Bor, 2.2} that
$X$ consists of regular elements of $S$ and is a right Ore set in both
$R$ and $S$. Consequently, the embedding of $R$ into $S$ extends to an
embedding of the right quotient ring $RX^{-1}$ into $SX^{-1}$. Note
that $RX^{-1}$ is artinian by Small's Theorem (see, e.g.,
\cite{\GooWar, 10.9}), that $SX^{-1}$ is finitely generated as a right
$RX^{-1}$-module, and that $SX^{-1}$ is the Goldie quotient ring of
$S$. Set $A = RX^{-1}$, $B = SX^{-1}$, and $K = QX^{-1}$.

As noted in (2.2iv), we can embed $A$ (as a left $A$-module) into a
finite direct sum of uniform left $A$-module factors of $A$ such that
the unique associated prime of each uniform factor is an associated
prime of $A$; let $L_1,\ldots,L_s$ denote these prime ideals. By
\cite{\LenLet, 2.3}, there now exist prime ideals $K_1,\ldots,K_t$
such that $0 = K_1K_2\cdots K_t A = K_1K_2\cdots K_t$ and such that
for each $1 \leq j \leq t$ there exists $1 \leq i \leq s$ for which
$L_i \idlink K_j$. However, since $K_1\cdots K_t \subseteq K$, and
since $K$ is a minimal prime ideal of $A$, it follows that $K = K_j$
for some $1 \leq j \leq t$. In particular, we have found an
annihilator prime $L$, of the left $A$-module $B$, such that $L
\idlink K$. Setting $Q' = L \cap R$, we see that $Q'$ is an
annihilator prime of $R$ for which there exists an ideal link from
$Q'$ to $Q$.

By (2.7i) and (4.1i), $\gkdim (R) = \gkdim (S) = \gkdim (R/Q) = \gkdim
(R/Q')$. Hence there exists a strong ideal link from $Q'$ to $Q$, by
(2.9ii). The lemma follows.  \qed\enddemo

The next lemma is another consequence of the Additivity Principle.

\proclaim{4.5 Lemma} Assume that $S$ is prime. 

{\rm (i)} Let $I$ be an ideal of $R$ not contained within any minimal prime
ideal of $R$. Then there exists a nonzero ideal of $S$ contained
within $SI$. If $S$ is finitely generated as a left $R$-module, then
there exists a nonzero ideal of $S$ contained within $IS$.

{\rm (ii)} Let $M$ be a faithful prime  $S$-module, and let $I$ be an
ideal of $R$ such that $IN = 0$ for some nonzero $R$-submodule $N$ of
$M$. Then $I$ is contained within a minimal prime ideal of
$R$. Consequently, if $Q$ is an annihilator prime of $M$ then $Q$ is a
minimal prime ideal of $R$. \endproclaim

\demo{Proof} (i) See \cite{\Letfour, 2.1}.

(ii) If $I$ is not contained within any minimal prime ideal of $R$
then, by (i), there exists a nonzero ideal $J$ of $S$ such that $J(SN)
= JN \subseteq SIN = 0$, a contradiction to the choice of $M$.
\qed\enddemo

The following is similar to \cite{\Lettwo, 2.5}.

\proclaim{4.6 Lemma} Assume that $S$ is prime and finitely generated
as both a left and right $R$-module. Let $M$ be a finitely generated,
uniform, faithful  $S$-module, and assume that $\Pa$ is the unique
associated prime of $M$. Then there exist prime ideals $\qa$ and $\qb$
of $R$ such that $\Pa$ lies over $\qa$, such that $\qb$ is a minimal
prime ideal of $R$, and such that $\qa \extlink \cdots \extlink
\qb$. \endproclaim

\demo{Proof} As observed in (2.2iv), $\Pa$ is the unique maximal
annihilator, in $S$, of nonzero $S$-submodules of $M$. Set $N = \ann
(\Pa)_M$. It follows from (4.5ii) that if $Q$ is an annihilator prime
of $N$, viewed as a left $R$-module, then $\Pa$ lies over $Q$. Next,
choose $L$ to be an $R$-submodule of $M$ maximal such that $L \cap N =
0$, and set $K = M/L$. The $R$-module embedding $N \hookrightarrow K$
is essential, forcing every associated prime, in $R$, of $K$ to also
be an associated prime of $N$. In particular, $\Pa$ lies over each
associated prime of $K$. Now let $X$ denote the union of the external
cliques of the associated primes of $K$. By (2.4), there exist
$Q_1,\ldots,Q_t \in X$ such that $Q_1\cdots Q_tK = 0$. Set $I =
Q_1\cdots Q_t$. Then $ISM = IM \subseteq L$.

Suppose that $I$ is not contained within a minimal prime ideal of $R$.
Then by (4.5i), there exists a nonzero ideal $J$ of $S$ contained
within $IS$. Note that $JM$ is nonzero because $M$ is faithful, and
that $JM$ is then an $S$-submodule of $M$ contained in $L$. Hence $JM
\cap N = 0$, contradicting the uniformity of $M$ as an
$S$-module. Thus $I$ is contained in a minimal prime ideal of $R$, and
one of the $Q_1,\ldots,Q_t$ must be a minimal prime ideal of $R$. The
lemma follows.  \qed\enddemo

\proclaim{4.7 Lemma} Suppose that $S$ has exact integer GK-dimension
and that $S$ is finitely generated as both a left and right
$R$-module. Let $M$ be a uniform  $S$-module. Let $\Pa$ be the
unique associated prime, in $S$, of $M$, and suppose that $\Pb$ is an
annihilator prime, in $S$, of $M$. Then there exist prime ideals $\qa$
and $\qb$ of $R$ such that $\qa$ lies under $\Pa$, such that $\qb$
lies under $\Pb$, and such that $\qa \extlink \cdots \extlink
\qb$. \endproclaim

\demo{Proof} First, we may assume without loss of generality that $\Pb
= \ann _SM$. Therefore, we can reduce to the case where $\Pb = 0$ and
$M$ is faithful as an $S$-module. By (2.6iv), we may further assume
that $M$ is finitely generated as a  $S$-module. The lemma now
follows from (4.6). \qed \enddemo

\proclaim{4.8 Proposition} Let $T$ be a noetherian $k$-algebra,
containing $S$ as a subalgebra, such that $T = t_1S + \cdots + t_\ell
S$.  Let $w$ be a positive integer. Assume there exist, for every
prime ideal $\qa$ of $R$, at most $w$ prime ideals $\qb$ such that
either $\qa$ is strongly ideal linked to $\qb$ or $\qb$ is strongly
ideal linked to $\qa$. Let $J$ be a prime ideal of $T$.

{\rm (i)} There exist at most $mw$ prime ideals of $S$ lying over
$Q$.

{\rm (ii)} There exist at most $mw$ prime ideals of $R$ lying
under $P$. 

{\rm (iii)} There exist at most $\ell m w$ prime ideals of $T$ lying
over $P$.

{\rm (iv)} There exist at most $\ell m^2 w$ prime ideals of $S$
lying under $J$. \endproclaim

\demo{Proof} Parts (i) and (ii) follow from (4.3) and (4.4).

(iii) Suppose that $J$ lies over $P$ and that $P$ lies over $Q$. Then
$J\cap R \subseteq Q$. It follows from (4.1i) that the GK-dimensions
of $T/J$, $S/P$, and $R/Q$ all coincide. Therefore, $Q$ must be
minimal over $J\cap R$, by (2.7ii), and so $J$ lies over $Q$. It
follows from (i) that at most $\ell m w$ prime ideals of $T$ can lie
over $Q$, and (iii) follows.

(iv) Suppose that $P$ lies under $J$ and directly over $Q$; as noted
in (4.1ii), $P$ lies directly over at least one prime ideal of $R$. As
in the proof of (iii), $Q$ lies under $J$. By (ii), there are no more
than $\ell m w$ prime ideals of $R$ lying under $J$. By (4.3i), there
are at most $m$ prime ideals of $S$ lying directly over $Q$. Part (iv)
follows.  \qed\enddemo

The following is similar, for example, to \cite{\Wartwo, 6.6}.

\proclaim{4.9 Lemma} Let $I/J$ be a strong ideal link from the prime
ideal $\qa$ to the prime ideal $\qb$ of $R$. Further suppose that $SJ
\cap R = J$ (e.g., suppose that $S$ is free as a right
$R$-module). Then there exist prime ideals $\Pa$ and $\Pb$ of $S$ such
that $\Pa$ lies over $\qa$, such that $\Pb$ lies over $\qb$, and such
there is a strong ideal link from $\Pa$ to $\Pb$. \endproclaim

\demo{Proof} Note that $\ann (S/SJ)_R = J$, and let $K = \ann
_S(S/SJ)$. By assumption, $K \cap R \subseteq J$, and so we may
suppose without loss of generality that $K = 0$. By (2.7i) and the
assumption that $I/J$ is a strong ideal link,
%%%%%%%%%%%%%%%%%%%%%%%%%%%%%%%%%%%%%%%%%%%%%%%%%%%%%%%%%%%%%%%%%%%%%%%
$$\gkdim (R) = \gkdim (S) = \gkdim (R/J) = \gkdim (R/\qa) = \gkdim
(R/\qb) .$$
%%%%%%%%%%%%%%%%%%%%%%%%%%%%%%%%%%%%%%%%%%%%%%%%%%%%%%%%%%%%%%%%%%%%%%%
Since $\gkdim (R/\qa) = \gkdim (R/\qb) = \gkdim (R)$, we see from
(2.7ii) that $\qa$ and $\qb$ are minimal prime ideals of $R$. It now
follows from \cite{\Wartwo, 6.3i} that there exist prime ideals $\Pa$
and $\Pb$ of $S$ such that $\Pa$ lies over $\qa$, such that $\Pb$ lies
over $\qb$, and such that there is an ideal link $L/M$ from $\Pa$ to
$\Pb$. Since $\gkdim (S/\Pa) = \gkdim (S/\Pb) = \gkdim(R/\qb) = \gkdim
(S)$, it follows from (2.9ii) that there is a strong ideal link from
$\Pa$ to $\Pb$. \qed\enddemo

\head 5. Extensions of Enveloping Algebras 
\endhead

In this section we establish bounds on the sizes of sets of
fundamental primes, in a setting that includes finite extensions of
enveloping algebras of reductive Lie algebras.

\subhead 5.1 \endsubhead Let $k$ be a field, and let $R$ be a
noetherian $k$-algebra of exact integer GK-dimension (see
(2.6)). Throughout this section, $S$ will denote a (necessarily
noetherian) $k$-algebra, containing $R$ as a subalgebra, such that $S$
is free of rank $\ell$ as a left $R$-module and is generated as a
right $R$-module by no more than $m$ elements. Further assume,
throughout, that every external clique (see (2.3ii)) of prime ideals
of $R$ has cardinality no greater than some fixed positive integer
$w$.

\subhead 5.2 \endsubhead For a general example fitting the description
of $R$ in (5.1), suppose that $k$ is an algebraically closed field of
characteristic zero, that $\a$ is a (finite dimensional) reductive
$k$-Lie algebra, and that $W$ is the associated Weyl group. Recall
that the enveloping algebra of a finite dimensional Lie algebra is
noetherian.

(i) Setting $U = U(\a)$, and $Z$ equal to the center of $U$, it
follows from \cite{\Soe, Theorem 1} that the sets of prime ideals of
$U$ having common intersection with $Z$ correspond exactly with
suitably defined $W$-orbits in the prime spectrum of the symmetric
algebra $S(\h ^\ast)$. In particular, by (2.5), the external cliques
in $\spec U$ have cardinality no greater than $\vert W \vert$.

(ii) It is shown, for example, in \cite{\KraLen, Chapter 6} that
enveloping algebras of finite dimensional Lie algebras, over any
field, have exact integer GK-dimension. So suppose, for the moment,
that $R$ is the enveloping algebra of a finite dimensional Lie
algebra. Recall, if $Q$ is a prime ideal of $R$, that $Q$ is left
primitive if and only if it is right primitive. It follows, when $P$
is a prime ideal of $S$, that $P$ is left primitive if and only if it
is right primitive (see, e.g., \cite{\Letone}).

(iii) Recall that the (left, or equivalently, right) primitive ideals
of $U$ (as defined in (i)) are precisely the prime ideals intersecting
to maximal ideals of $Z$, that the minimal primitive ideals of $U$ are
the ideals of the form $U\chi$ for maximal ideals $\chi$ of $Z$, and
that the set of annihilators of Verma modules over $U$ coincides with
the set of minimal primitive ideals of $U$. The Verma modules are
uniform $U$-modules of finite length, and from Duflo's Theorem it
follows that every primitive ideal is the annihilator of a simple
factor of a Verma module. (See \cite{\Dix} or \cite{\Jan} for
details.) It now follows from (2.4) that the external clique in $\spec
U$ of a primitive ideal $P$ of $U$ is exactly the set of primitive
ideals whose intersection with $Z$ equals $P\cap Z$.

(iv) It follows from (iii) that if $\qa$ and $\qb$ are primitive
ideals of $U$ such that $\qa \subseteq \qb$, then $\qa$ and $\qb$ are
in the same external clique. Suppose that $M$ is a uniform $U$-module
with associated prime $Q$, that $Q$ is primitive, and that $L$ is a
simple $U$-module factor of $M$. It now follows from (2.12iii) that
the annihilator in $U$ of $L$ is in the external clique of $Q$.

\subhead 5.3 \endsubhead Let $T = M_\ell (R)$. By \cite{\Lor, 1.6},
both $S$ and $T$ have exact integer GK-dimension. In
particular, any pair of strongly ideal linked prime ideals of $R$,
$S$, or $T$ is contained within a common clique, by (2.9iii). As noted
in \cite{\Lettwo, 7.1}, there is a $k$-algebra embedding of $S$ into
$T$ such that $T$ is a free right $S$-module of rank $\ell$. Of
course, every external clique of prime ideals of $T$ has cardinality
no greater than $w$.

\proclaim{5.4 Lemma} Let $P$ be a prime ideal of $S$.

{\rm (i)} Let $C$ denote the union of the external cliques of the
prime ideals of $R$ lying under $P$, and let $C'$ denote the set of
prime ideals in $S$ lying over ideals in $C$. Then $\vert C' \vert
\leq m^2 w$.

{\rm (ii)} Let $D$ denote the union of the external cliques of the
prime ideals of $T$ lying over $P$, and let $D'$ denote the set of
prime ideals in $S$ lying under ideals in $D$. Then $\vert D' \vert
\leq \ell^2m^3 w^3$. \endproclaim

\demo{Proof} (i) By (4.3ii) and (4.4), $\vert C \vert \leq mw$. By
(4.4), $C'$ is equal to the set of prime ideals of $S$ lying directly
over ideals in $C$, and so from (4.3i) it follows that $\vert C' \vert
\leq m^2w$.

(ii) By (4.8iii), $\vert D \vert \leq \ell m w^2$, and so it follows from
(4.8iv) that $\vert D' \vert \leq \ell^2 m^3 w^3$. \qed\enddemo

Fundamental series and primes (cf\. \cite{\BroWar, \Jat, \LenWar}) are
discussed in (2.6) and (2.12). 

\proclaim{5.5 Theorem} If $P$ is a prime ideal of $S$ then $\vert
\Fund (P) \vert \leq \ell^2 m^5 w ^4$. \endproclaim

\demo{Proof} First, let $A$ denote the set of annihilator primes of
uniform  $S$-modules with associated prime $P$. Set $C$ equal to
the union of the external cliques of prime ideals of $R$ lying under
$P$, and let $C'$ denote the set of prime ideals of $S$ lying over
ideals in $C$. By (4.7), $A \subseteq C'$, and so $\vert A \vert \leq
m^2w$, by (5.4i).

Second, let $\Pa$ be a prime ideal of $S$, and let $B$ denote the set
of prime ideals $\Pb$ of $S$ for which there is a strong ideal link
from $\Pa$ to $\Pb$. Let $D$ be the union of the external cliques of
the prime ideals of $T$ lying over $\Pa$, and let $D'$ denote the set
of prime ideals of $S$ lying under prime ideals in $D$.  By (2.12i)
and (4.9), $B \subseteq D'$. Thus $\vert B \vert \leq \ell^2 m^3 w^3$,
by (5.4ii).

It follows from (2.8) that $\vert \Fund (P) \vert \leq \vert A \vert \
\vert B \vert$. The theorem follows. \qed\enddemo

\subhead 5.6 \endsubhead A bound similar to, and sharper than, that of
(5.5) can be established for external links (see (2.3)) and can in
turn be applied to extensions of simple modules, as follows: Suppose
that $0 \rightarrow \La \rightarrow M \rightarrow \Lb \rightarrow 0$
is a nonsplit extension of simple modules, that $\Pa = \ann _S\La$,
and that $\Pb = \ann _S\Lb$; then $M$ is an external link from $\Pa$
to $\Pb$.

\proclaim{5.7 Theorem} Let $P$ be a prime ideal of $S$. 

{\rm (i)} There exist at most $\ell^2 m^3 w^3$ prime ideals $P'$ of
$S$ such that $P \link P'$.

{\rm (ii)} There exist at most $m^2w$ prime ideals $P'$ of $S$ such
that $P \extlink P'$ but $P \notlink P'$.

{\rm (iii)} There exist at most $m^2 w + \ell^2 m^3 w^3$ prime ideals
$P'$ of $S$ such that $P \extlink P'$. \endproclaim

\demo{Proof} (i) By (2.10) and (4.9), if $P \link P'$ then there exist
prime ideals $J$ and $J'$ of $T$ such that $J$ lies over $P$, such
that $J'$ lies over $P'$, and such that there is a strong ideal link
from $J$ to $J'$. Part (i) now follows from (2.12i) and (5.4ii).

(ii) Let $M$ be an external link from $P$ to $P'$, and suppose that $P
\notlink P'$. By (2.3iii), $P' = \ann _SM$, and so by (4.6), or (4.7),
there exist prime ideals $Q$ and $Q'$ of $R$ such that $P$ lies over
$Q$, $P'$ lies over $Q'$, and $Q \extlink \cdots \extlink Q'$. Part
(ii) now follows from (5.4i).

(iii) This part now follows from (i), (ii), and (2.3iii). \qed\enddemo

\subhead 5.8 \endsubhead Assume that $R = U(\a )$, for $\a$ as in (5.2).

(i) Distinct primitive ideals with finite codimension in $R$ cannot
have equal intersection with its center, and the analogous statement
also holds true for $T$. Hence, there are no ideal links in either $R$
or $T$ between distinct finite codimensional primitive ideals. It now
follows, for example, from (4.9) that linked, primitive, finite
codimensional ideals in $S$ lie under a common finite codimensional
primitive ideal of $T$. Moreover, there is a nonsplit extension of
finite dimensional simple $S$-modules if and only if there exists a
direct link between their corresponding annihilators.

(ii) Suppose that $J$ is a finite codimensional primitive ideal of
$T$, that $P$ is a finite codimensional primitive ideal of $S$, and
that $Q$ is finite codimensional primitive ideal of $R$. By (i) and
(4.4), $Q$ lies under $J$ (resp\. $P$) if and only if $Q$ lies
directly under $J$ (resp\. $P$). Therefore, by (4.3ii) there exist at
most $\ell m$ prime ideals of $R$ lying under $J$, and by (4.3i) there
exist at most $m$ prime ideals of $S$ lying over $Q$. Thus, as in the
proof of (4.8), there exist at most $\ell m^2$ prime ideals of $S$ lying
under $J$.

(iii) Let $L$ be an arbitrary finite dimensional simple 
$S$-module, and let $X$ denote the set of finite dimensional simple
 $S$-modules $L'$ for which either $\Ext _S(L, L') \ne 0$ or
$\Ext _S (L', L) \ne 0$. It follows from (i) and (ii) that $\vert X
\vert \leq \ell m^2$.

(iv) Let $\Pa$ and $\Pb$ be prime ideals of $S$, and suppose that $M$
is an external link from $\Pa$ to $\Pb$. We can show as follows that
$\Pa$ is primitive if and only if $\Pb$ is primitive. First, if $\Pa
\link \Pb$, then $\Pa$ is primitive if and only if $\Pb$ is primitive,
by (2.3i) and (5.2ii). We may therefore assume that $\Pa \notlink
\Pb$. By (2.3iii) $\Pb = \ann_SM$, and by (4.7) it follows that there
exist prime ideals $\qa$ and $\qb$ of $R$ such that $\Pa$ lies over
$\qa$, such that $\Pb$ lies over $\qb$, and such that $\qa \extlink
\cdots \extlink \qb$. It now follows from (5.2iii) and (4.1i) that
$\Pa$ is primitive if and only if $\Pb$ is primitive.

\subhead 5.9 \endsubhead Retaining the assumptions of (5.1), suppose
further that the characteristic of $k$ is not equal to two, that $S$
is ($\Z_2$-)graded, that $R=\Ro$ is a trivially graded subalgebra of
$\So$, and that the free left $R$-module basis for $S$ is comprised of
homogeneous elements. The algebra $T$ is then also graded, and $R$ and
$S$ are graded subalgebras of $T$. Following (3.6), the results in
this section remain valid when the modules and bimodules involved are
replaced with their graded counterparts.

\head 6. Classical Lie Superalgebras \endhead

We now apply the results of the preceding sections to the
representation theory of Lie superalgebras. The reader is referred to
\cite{\Kac} or \cite{\Sch} for background information.

\subhead 6.1 \endsubhead (i) Let $k$ be an algebraically closed field
of characteristic zero. Assume throughout this section that $\g =
\goi$ is a finite dimensional $k$-Lie superalgebra, that $\V$ is the
associative enveloping algebra of $\g$, and that $\U$ is the
enveloping algebra of the finite dimensional Lie algebra $\go$.

(ii) It follows from the Poincar\/e-Birkhoff-Witt Theorem for Lie
superalgebras that $\V$ contains $\U$ as a subalgebra and that $\V$ is
free of rank $2^{\dim\gi}$ as both a right and left $\U$-module (and
that the right and left bases can be chosen to coincide). Hence $\V$
is noetherian. As remarked in (5.3), $\U$ and $\V$ have exact integer
GK-dimension.

(iii) The ($\Z_2$-)grading on $\g$ extends to a grading of $\V$, and
$\U \subseteq \Vo$. Furthermore, the free left and right common
$\U$-basis for $\V$ can be chosen to consist of homogeneous
elements. In particular, the remarks in (5.9) apply.

\subhead 6.2 \endsubhead (i) It follows from (5.2ii) that a prime
ideal of $\V$ is left primitive if and only if it is right primitive.

(ii) Assume that $\g$ is simple (i.e., it contains no nonzero graded
ideals other than itself). If $\gi$ is completely reducible as an
$\ad\go$-module -- a condition that occurs if and only if $\go$ is
reductive (see, e.g., \cite{Sch, Theorem 1, p\. 101}) -- then $\g$ is
called {\sl classical\/}.

(iii) Suppose that $\go$ is reductive (with $\g$ not necessarily
simple), and let $P$ be a (left or right) primitive ideal of $\V$. It
follows from (5.8iv) that the external clique, in $\spec \V$, of $P$
consists only of primitive ideals, and so it follows from (2.12i) that
the prime ideals in $\Fund (P)$ are all primitive. Similarly, if $p$
is a graded-primitive ideal of $\V$, then the graded-prime ideals in
$\Grfund (p)$ are graded-primitive.

\subhead\nofrills\endsubhead 

The next theorem summarizes the immediate consequences of (5.2),
(5.5), (5.7iii), and (5.8iii). Passage to the graded case follows from
(3.6) and (5.9).

\proclaim{6.3 Theorem} Assume that $\go$ is reductive (e.g., $\g$ is
classical simple) and that the Weyl group associated to $\go$ has
cardinality $w$. Set $\ell = 2^{\dim\gi}$. Let $P$ be a prime
ideal of $\V$, and let $p$ be a graded-prime ideal of $\V$.

{\rm (i)} $\vert \Fund (P) \vert \leq \ell^7 w^4$.  

{\rm (i$'$)} $\vert \Grfund (p) \vert \leq \ell^7 w^4$.

{\rm (ii)} There exist at most $\ell^2 w + \ell^5 w^3$ prime ideals
$P'$ of $\V$ such that $P \extlink P'$. 

{\rm (ii$'$)} There exist at most $\ell^2 w + \ell^5 w^3$ graded-prime ideals
$p'$ of $\V$ such that $p \extlink p'$. 

{\rm (iii)} Let $L$ be a simple finite dimensional $\V$-module. Then
there are at most $\ell^3$ finite dimensional simple $\V$-modules $L'$
for which $\Ext (L,L')$ or $\Ext (L',L)$, calculated in $\mod \V$, is
not equal to zero.

{\rm (iii$'$)} Let $L$ be a graded-simple finite dimensional
$\V$-module. Then there are at most $\ell^3$ finite dimensional
graded-simple $\V$-modules $L'$ for which $\Ext (L,L')$ or $\Ext
(L',L)$, calculated in $\grmod \V$, is not equal to zero.
\endproclaim

\subhead 6.4 \endsubhead We now turn to highest weight modules,
following \cite{\Musone} (cf\. \cite{\PenSer}). Assume that $\g =
\goi$ is classical simple, that $W$ is the Weyl group of $\go$, that
$w = |W|$, and that $\ell = 2^{\dim\gi}$.

(i) To start, fix a triangular decomposition $\g = \nm \oplus \h
\oplus \np$ as in \cite{\Musone, 1.1}. Then $\go = \nmo \oplus \ho
\oplus \npo$ is a triangular decomposition of $\go$.

(ii) Set $\b = \h \oplus \np$. For each $\lambda \in \ho ^\ast$ there
is a unique finite dimensional graded-simple $\b$-module
$K_\lambda$ such that $\np . K_\lambda = 0$, and such that $h.x =
\lambda(h).x$ for all $h \in \ho$ and $x \in K_\lambda$. Moreover, all
of the finite dimensional graded-simple $\b$-modules are obtained, up
to isomorphism, in this fashion. (See \cite{\Kac, \S 5.2},
\cite{\Musone, \S 1.1}.)

(iii) Fix $\lambda \in \ho ^\ast$, and set $\M (\lambda) = \V \otimes
_{U(\b )}K_\lambda$. Then $\M (\lambda )$ has finite length as a
$\V$-module and has a unique maximal graded $\V$-submodule; the
resulting graded-simple factor module is denoted $\L (\lambda)$.

(iv) Extending Duflo's Theorem (see, e.g., \cite{\Jan, 7.4}), Musson
proved that the map $\ho^\ast \rightarrow \grprim\V$, sending $\lambda
\mapsto \ann _\V \L (\lambda)$, is surjective \cite{\Musone, \S
2}. Moreover, the fibers of the preceding map all have cardinality no
greater than $\ell w$; see \cite{\Letfour, 3.2}

\subhead\nofrills\endsubhead 

In view of (5.6), the next proposition follows immediately from
(6.3ii$'$) and (6.4iv).

\proclaim{6.5 Proposition} Retain the notation and assumptions of {\rm
(6.4)}. Let $\lambda \in \ho^\ast$. Then there are at most $\ell^3 w^2
+ \ell^6 w^4$ linear forms $\lambda' \in \ho^\ast$ for which
%%%%%%%%%%%%%%%%%%%%%%%%%%%%%%%%%%%%%%%%%%%%%%%%%%%%%%%%%%%%%%%%%%%
$$\Ext \left(\L(\lambda),\L(\lambda')\right) \quad \text{or} \quad
\Ext \left(\L(\lambda'),\L(\lambda)\right),$$
%%%%%%%%%%%%%%%%%%%%%%%%%%%%%%%%%%%%%%%%%%%%%%%%%%%%%%%%%%%%%%%%%%%
calculated in $\grmod \V$, is not equal to zero. \endproclaim

We conclude our study by considering graded-essential extensions of
graded-simple highest weight modules. 

\proclaim{6.6 Lemma} Retain the notation and assumptions of {\rm
(6.4)}, and let $P$ be a primitive ideal of $\V$. Then there exist at
most $\ell^2 w^3$ primitive ideals $P'$ of $\V$ such that $P \subseteq
P'$. If $p$ is a graded-primitive ideal of $\V$ then there exist at
most $\ell^2 w^3$ graded-primitive ideals $p'$ of $\V$ such that $p
\subseteq p'$. \endproclaim

\demo{Proof} Suppose that $Q_1,\ldots,Q_t$ are the primitive (see
(4.1i)) ideals of $\U$ lying under $P$. By (4.8ii), $t \leq \ell
w$. Now suppose that $P'$ is a primitive ideal of $\V$ containing $P$,
and suppose that $P$ lies over the primitive ideal $Q'$ of $\U$. Then
$Q'$ contains $Q_i$, for some $1 \leq i \leq t$, since there will be a
power of $Q_1\cdots Q_t$ contained inside $P\cap \U \subseteq P' \cap
\U \subseteq Q'$. However, there are at most $w$ primitive ideals of
$\U$ containing $Q_i$, and at most $\ell w$ prime ideals of $\V$ lying
over $Q'$, by (4.8i). The ungraded part of the lemma now follows, and
the graded case is proved similarly. \qed\enddemo

\proclaim{6.7 Theorem} Retain the notation and assumptions of {\rm
(6.4)}. Let $M$ be a graded-uniform $\V$-module whose unique
associated graded-prime ideal is graded-primitive (e.g., let $M$ be a
graded-indecomposable graded-injective $\V$-module with a nonzero
graded-socle). Then there exist at most $\ell^{10} w^8$ linear forms
$\lambda \in \ho^\ast$ such that $\L (\lambda)$ is isomorphic to a
graded-simple $\V$-module subfactor of $M$. \endproclaim

\demo{Proof} Suppose that $\L (\lambda)$ is a graded-simple
$\g$-subfactor of $M$, for some $\lambda \in \ho^\ast$, and set
$\plambda = \ann_\V \L (\lambda)$. By (2.12iii), $\plambda$ contains
some graded-prime ideal $p \in G = \Grfund (M)$. By (6.3i$'$), $\vert
G \vert \leq \ell^7 w^4$. But $p$ is graded-primitive, by (6.2iii),
and so there exist at most $\ell^2 w^3$ graded-primitive ideals of
$\V$ containing $p$, by (6.6). However, as noted in (6.4iv), there
exist at most $\ell w$ graded-simple heighest weight modules with a
given graded-primitive annihilator in $\V$. The theorem
follows. \qed\enddemo

\remark{6.8 Remark} (i) I do not know the sharpness, in general, of the
various bounds established in this section. Of course, when $\g = \go$
they significantly overestimate the situation. \endremark

(ii) Musson's example \cite{\Mustwo, \S 4}, described in (1.2), shows
that $\V$ may have infinite cliques (see (2.3)) of finite
codimensional graded-primitive ideals.

\enddocument